\documentstyle[11pt]{amsart}
\newtheorem{th}{Theorem}

\newtheorem{prp}{Proposition}
\newtheorem{lm}{Lemma}
\newtheorem{rmk}{Remark}

\numberwithin{equation}{section}

\begin{document}
\title[Cohomology of simple modules for
${\mathfrak{sl}}_3(k)$]{Cohomology of simple modules for ${\mathfrak{sl}}_3(k)$ in
characteristic $3$}
\author{A.A.Ibrayeva$^1$, Sh.Sh.Ibraev$^2$, G.K.Yeshmurat$^3$}
\address{$^{1,2,3}$Korkyt Ata Kyzylorda University, Kyzylorda, Kazakhstan}
\email{ibrayevsheraly@@gmail.com}
\date{}
\begin{abstract} In this paper, we calculate cohomology of a classical Lie algebra of
type $A_2$ over an algebraically field $k$ of characteristic $p=3$
with coefficients in simple modules. To describe their structure, we
will consider them as modules over an algebraic group $SL_3(k).$ In
the case of characteristic $p=3,$ there are only two peculiar simple
modules: a simple module isomorphic to the quotient module of the
adjoint module by the center, and a one-dimensional trivial module.
The results on the cohomology of simple nontrivial module are used
for calculate the coomology of the adjoint module. We also calculate
cohomology of the simple quotient algebra Lie of $A_2$ by the
center.
\vspace{0,1cm}


\vspace{0,1cm}


\end{abstract}
\maketitle \vspace{0,1cm}

\vspace{0,1cm}

\section{Introduction}\label{s1}

The cohomology theory of modular Lie algebras is one of the
interesting questions in the theory of Lie algebras. Many
significant results are devoted to the study of the cohomology of
classical modular Lie algebras. Their restricted cohomology with
coefficients in the dual Weyl modules was studied in \cite{1} -- \cite{3}.
Central extensions are described in \cite{4}, \cite{5}. In \cite{6} and \cite{7} the
outer derivations are calculated. As the second cohomology, local
deformations are calculated in \cite{8}--\cite{10}.

 Among the classical modular
Lie algebras, the cohomology of simple modules is completely
described only for a three-dimensional Lie algebra of type $A_1$
\cite{11}. It is known that, for other classical modular Lie algebras a
complete description of the cohomology of simple modules has not yet
been obtained. In this paper, we give a complete description of such
cohomology for the Lie algebra of type $A_2$ over an algebraically
closed field of characteristic $p=3.$ The first cohomology groups of
simple modules for $A_2$ was computed in \cite{12}. A similar result for
the second cohomology groups was obtained in \cite{13}. In all other
cases the computation of the cohomology structure of simple modules
for $A_2$ is close to completion. The results will be published in
the next works of the second author.

Let us introduce the basic definitions and notation. Let
$\mathfrak{g}$ be a Lie algebra over a field $k$ characteristics of
$p$ and $ M $ be a $\mathfrak{g}$-module. We denote the $n$-th
exterior power of the space $\mathfrak{g}$ by
$\Lambda^n({\mathfrak{g}})$ and let
$$C^n({\mathfrak{g}},M)=Hom(\Lambda^n,M)=\langle\;\psi\,:\;{\mathfrak{g}}\times
\cdots\times{\mathfrak{g}}\to M\;\rangle_k,\;n>0$$ is a space of
multilinear skew-symmetric mappings in $n$ arguments with
coefficients in $M.$ We put
$$C^n({\mathfrak{g}},M)=0,\;n<0,\;\;C^0({\mathfrak{g}},M)=M,\;\;C^*({\mathfrak{g}},M)=
\bigoplus_{n=-\infty}^{+\infty}C^n({\mathfrak{g}},M).$$ Define the
coboundary operator $$\,d\,:\;C^*({\mathfrak{g}},M)\longrightarrow
C^*({\mathfrak{g}},M)$$ as follows
$$d\psi(l_1,l_2,\cdots
,l_{n+1})=$$
$$
\sum_{i<j}(-1)^{i+j}\psi([l_i,l_j],\cdots , \hat{l_i},\cdots
,\hat{l_j},\cdots, l_{n+1})+$$$$
\sum_i(-1)^{i+1}[l_i,\psi(l_1,\cdots,\hat{l_i},\cdots,l_{n+1})],$$
where $\psi\in C^n({\mathfrak{g}},M).$ Then $d^2=0,$ therefore
$B^*({\mathfrak{g}},M)\subseteq Z^*({\mathfrak{g}},M),$ where
$$Z^*({\mathfrak{g}},M)=\langle\,\psi\in
C^*({\mathfrak{g}},M)\;:\;d\psi=0\;\rangle_k,$$
$$B^*({\mathfrak{g}},M)=\langle\,d\psi\;:\;\psi\in C^*({\mathfrak{g}},M)\;\rangle_k.$$
So, we can introduce the factor-space
$$H^*({\mathfrak{g}},M)=Z^*({\mathfrak{g}},M)/B^*({\mathfrak{g}},M).$$
The spaces
$C^*({\mathfrak{g}},M),\;Z^*({\mathfrak{g}},M),\;B^*({\mathfrak{g}},M),\;
H^*({\mathfrak{g}},M)$  are called {\it space of cochains}, {\it space
of cocycles}, {\it spaces of coboundaries}, and {\it space of
cohomologies} of the Lie algebra ${\mathfrak{g}}$ with coefficients in
the ${\mathfrak{g}}$-module $M$ respectively.

Similarly, the spaces
$$C^n({\mathfrak{g}},M),\;Z^n({\mathfrak{g}},M)=Z^*({\mathfrak{g}},M)\cap
C^n({\mathfrak{g}},M),\;$$
$$B^n({\mathfrak{g}},M)=B^*({\mathfrak{g}},M)\cap
C^n({\mathfrak{g}},M)\;\mbox{ט}\;H^n({\mathfrak{g}},M)=H^*({\mathfrak{g}},M)\cap
C^n({\mathfrak{g}},M)$$ are called {\it space of $n$-cochains}, {\it
space of $n$-cocycles}, {\it spaces of $n$-coboundaries}, and {\it
space of $n$-cohomologies} of the Lie algebra ${\mathfrak{g}}$ with
coefficients in the ${\mathfrak{g}}$-module $M$ respectively.

We call the ${\mathfrak{g}}$-module $M$ is {\it peculiar}, if
$H^*({\mathfrak{g}},M)\neq 0.$ We say that $M$ is {\it $n$-peculiar}
module over ${\mathfrak{g}},$ if $H^n({\mathfrak{g}},M)\neq 0.$

Now let ${\mathfrak{g}}$ be a classical Lie algebra of type $A_2$ over
algebraically closed field $k$ of positive characteristic $p>0$ and
$M$ is a ${\mathfrak{g}}$-module. We decompose $C^{*}({\mathfrak{g}},M)$
into a direct sum of weight subspaces with respect to the maximal
torus $T$ of the group $G=SL_3(k):$
$$C^{*}({\mathfrak{g}},M)=\bigoplus_{\mu\in X(T)}C^{*}_{\mu}({\mathfrak{g}},M),$$
where $X(T)$ is the additive character group of $T.$ Then
$$H^n({\mathfrak{g}},M)=\bigoplus_{\mu\in X(T)}H^n_{\mu}({\mathfrak{g}} M).$$

Identify the space $C^n({\mathfrak{g}},M)$  with the space
$\bigwedge^n{\mathfrak{g}}^*\bigotimes M$ and denote by $\prod(V)$ the
set of weights of the $G$-module subspace $V$ of
$H^{*}({\mathfrak{g}},M).$

Since $\prod(H^n({\mathfrak{g}},M))\subseteq
pX(T)\bigcap\prod(\bigwedge^n{\mathfrak{g}}^*\bigotimes M),$ then we
can consider only the elements of the subspace
$\overline{C}^n({\mathfrak{g}},M)$ of $C^n({\mathfrak{g}},M)$ with
weights contained in the set
$pX(T)\bigcap\prod(\bigwedge^n{\mathfrak{g}}^*\bigotimes M).$ The
corresponding subspaces of cocycles and cohomologies are denoted by
$\overline{Z}^n({\mathfrak{g}},M)$  and
$\overline{H}^n({\mathfrak{g}},M).$  Note that
$$H^n({\mathfrak{g}},M)=\overline{H}^n({\mathfrak{g}},M).$$ We will use the following
well known formulas:
$$\dim\,H^n({\mathfrak{g}},M)=\dim\,\overline{Z}^n({\mathfrak{g}},M)+\dim\,
\overline{Z}^{n-1}({\mathfrak{g}},M)-\dim\,\overline{C}^{n-1}({\mathfrak{g}},M),
\eqno(1)$$
$$\dim\,H^n({\mathfrak{g}},M)=\dim\,H^{\dim\,{\mathfrak{g}}-n}({\mathfrak{g}},M^*).
\eqno(2)$$ The weight subspaces are invariant under the action of
the coboundary operator, therefore the formula~(1)is also holds for
weight subspaces:
$$
\dim\,H^n_{\mu}({\mathfrak{g}},M)=\dim\,\overline{Z}^n_{\mu}({\mathfrak{g}},M)+\dim\,
\overline{Z}^{n-1}_{\mu}({\mathfrak{g}},M)-\dim\,\overline{C}^{n-1}_{\mu}({\mathfrak{g}},M).
\eqno(3)$$

Let $L(r,s)$ denote a simple ${\mathfrak{g}}$-module with the highest
weight $r\omega_1+s\omega_2,$ where $\omega_1,\,\omega_2$ are
fundamental weights.

It is known that the composition of a representation of $SL_3(k)$ on
a vector space $L$ with a $d$-th power of the Frobenius map defines
a new representation on which the Lie algebra ${\mathfrak{g}}$ acts
trivially. We denote the resulting module by $L^{(d)}.$ To each
weight $\mu$ of the space $L$ there corresponds a weight $p^d\mu$ of
the space $L^{(d)}.$ The cohomology group $H^n({\mathfrak{g}},M),$ as
a $SL_3(k)$-module, consists either a twisted module $L^{(d)}$ for
some $d,$ or a one-dimensional trivial module $k.$ For the
multiplicity of a $SL_3(k)$-module $L^{(d)}$ in
$H^n({\mathfrak{g}},M)$, we will use the notation
$[H^n({\mathfrak{g}},M):L^{(d)}].$ Further, for convenience, we will
use the following abbreviations:
$H^n({\mathfrak{g}},k):=H^n({\mathfrak{g}}),$ $\bigoplus_{i=1}^mV:=mV,$
where $V$ is a $SL_3(k)$-module.

Let's formulate the main result of this paper:

\begin{th}\label{th1} Let ${\mathfrak{g}}$ be a classical Lie algebra of
type $A_2$ over an algebraically closed field $k$ of characteristic
$p=3$ and $M$ ba a simple ${\mathfrak{g}}$-module. Then there are the
following isomorphisms of $SL_3(k)$-modules:

(a) $H^0({\mathfrak{g}})\cong H^8({\mathfrak{g}})\cong k,$
$H^2({\mathfrak{g}})\cong H^6({\mathfrak{g}})\cong L(1,0)^{(1)}\oplus
L(0,1)^{(1)},$ $H^3({\mathfrak{g}})\cong H^5({\mathfrak{g}})\cong
L(1,0)^{(1)}\oplus L(0,1)^{(1)}\oplus k;$

(b) $H^1({\mathfrak{g}},L(1,1))\cong H^7({\mathfrak{g}},L(1,1))\cong
L(1,0)^{(1)}\oplus L(0,1)^{(1)}\oplus k,$

$H^3({\mathfrak{g}},L(1,1))\cong H^5({\mathfrak{g}},L(1,1))\cong
H^0(1,1)^{(1)},$ $H^4({\mathfrak{g}},L(1,1))\cong 2H^0(1,1)^{(1)}.$

In other cases $H^n({\mathfrak{g}},M)=0.$
\end{th}

\section{Proof of the Theorem~\ref{th1}}\label{s2} As the basis vectors for ${\mathfrak{g}}$ we choose the special
derivations of the algebra of divided powers $O_3(\mathbf{1}):$
$$h_1=x_1\partial_1-x_2\partial_2,
h_2=x_2\partial_2-x_3\partial_3, e_1=x_1\partial_2,
e_2=x_2\partial_3, e_3=x_1\partial_3,$$$$ f_1=x_2\partial_1,
f_2=x_3\partial_2, f_3=x_3\partial_1.$$ Over a field of
characteristic $p=3,$ the Lie algebra ${\mathfrak{g}}$ is not simple,
it has a one-dimensional center $\langle h_1-h_2\rangle_k.$ The
quotient algebra by the center is a simple Lie algebra; we denote it
by $\overline{{\mathfrak{g}}}$ or $\overline{A_2}.$

It is known that the peculiar modules of the Lie algebra
${\mathfrak{g}}$ are restricted \cite{11}. According to Lemma~3.1 in \cite{13},
only the following two simple restricted modules are peculiar:
$L(0,0)\cong k$ and $L(1,1)\cong \overline{{\mathfrak{g}}}.$ For
$L(1,1)$ we get the following description:
$$L(1,1)\cong \langle h_1,\,h_2, e_1,\,e_2,\,e_3,\,f_1,\,f_2,\,f_3\,:\,h_1-h_2=0\rangle_{k}.$$ Consider each of these modules separately.

Let $M=L(0,0)\cong k.$

\begin{lm}\label{l1} There are the following isomorphisms of
$SL_3(k)$-modules:

(a) $H^0({\mathfrak{g}})\cong k;$

(b) $H^2({\mathfrak{g}})\cong L(1,0)^{(1)}\oplus L(0,1)^{(1)};$

(c) $H^3({\mathfrak{g}})\cong L(1,0)^{(1)}\oplus L(0,1)^{(1)}\oplus
k;$

(d) $H^5({\mathfrak{g}})\cong L(1,0)^{(1)}\oplus L(0,1)^{(1)}\oplus
k;$

(e) $H^6({\mathfrak{g}})\cong L(1,0)^{(1)}\oplus L(0,1)^{(1)};$

(f) $H^8({\mathfrak{g}})\cong k.$

In other cases $H^n({\mathfrak{g}})=0.$
\end{lm}

{\bf Proof.} The statements $(a)$ and $(f)$ are obvious. The
triviality of $H^1({\mathfrak{g}})$  in characteristic $p=3$ was
proved in \cite{12}.

$(b)$ The set $\prod(\overline{C}^2({\mathfrak{g}}))$ consists only
the following weights:
$$0\,\pm3\omega_1,\,\pm3(\omega_1-\omega_2),\,\pm3\omega_2.$$
Therefore, only the trivial one-dimensional module and the twisted
simple modules $L(1,0)^{(1)},$ $L(0,1)^{(1)},$ can be as nonzero
composition factors of $H^2({\mathfrak{g}}).$ They are generated by
the classes of cocycles with dominant weights $0,$ $3\omega_1,$ and
$3\omega_2$ respectively.

The subspace $\overline{C}_{0}^2({\mathfrak{g}})$ is $4$-dimensional
and spans  by the cochains $$h_1^*\wedge h_2^*,\,e_1^*\wedge f_1^*,\,
e_2^*\wedge f_2^*,\, e_3^*\wedge f_3^*.$$ If $$a_1h_1^*\wedge
h_2^*+a_2e_1^*\wedge f_1^*+a_3e_2^*\wedge f_2^*+a_4e_3^*\wedge
f_3^*\in \overline{Z}^2({\mathfrak{g}})$$ then, by cocycle condition,
$a_1=0,\,a_4=a_2+a_3.$ Therefore
$\dim\,\overline{Z}_{0}^2({\mathfrak{g}})=2.$ Since
$\dim\,\overline{C}_{0}^1({\mathfrak{g}})=2$ and
$\dim\,\overline{Z}_{0}^1({\mathfrak{g}})=0,$ by~(2),
$$\dim\,\overline{H}_{0}^2({\mathfrak{g}})=2+0-2=0.$$

The subspace $\overline{C}_{3\omega_1}^2({\mathfrak{g}})$ is
one-dimensional and spans  by the cochain $f_1^*\wedge f_3^*.$
Notice that $af_1^*\wedge f_2^*\in\overline{Z}^2({\mathfrak{g}})$ for
all $a\in k.$ Therefore
$\dim\,\overline{Z}_{3\omega_1}^2({\mathfrak{g}})=1.$ Since
$\dim\,\overline{C}_{3\omega_1}^1({\mathfrak{g}})=0,$  by~(3),
$\dim\,\overline{H}_{3\omega_1}^2({\mathfrak{g}})=1.$ So,
$[H^2({\mathfrak{g}}):L(1,0)^{(1)}]=1.$

Arguing as in the previous case, we obtain
$[H^2({\mathfrak{g}}):L(0,1)^{(1)}]=1.$ Thus $H^2({\mathfrak{g}})\cong
L(1,0)^{(1)}\oplus L(0,1)^{(1)}.$

$(c)$ The sets of weights $\prod(\overline{C}^3({\mathfrak{g}}))$ and
$\prod(\overline{C}^2({\mathfrak{g}}))$ are coincide. Therefore, we
consider only the weight subspaces of $3$-cochains corresponding to
the dominant weights $0,$ $3\omega_1,$ and $3\omega_2.$

The subspace $\overline{C}^3_{0}({\mathfrak{g}})$ is $8$-dimensional
and spans  by the cochains
$$h_1^*\wedge e_1^*\wedge f_1^*,\,h_2^*\wedge e_1^*\wedge f_1^*,\,
h_1^*\wedge e_2^*\wedge f_2^*,\, h_2^*\wedge e_2^*\wedge f_2^*,$$
$$h_1^*\wedge e_3^*\wedge f_3^*,\, h_2^*\wedge e_3^*\wedge f_3^*,\,
e_3^*\wedge f_1^*\wedge f_2^*,\, e_1^*\wedge e_2^*\wedge f_3^*.$$
Suppose that a linear combination of these vectors with coefficients
$b_i,\,i= 1,\cdots,8$  respectively, is a $3$-cocycle. Then the
cocycle condition implies that

$b_1+b_2+b_5+b_7-b_8=0,$

$b_2+b_3-b_7+b_8=0,$

$b_3+b_4+b_6+b_7-b_8=0,$

$2b_4+2b_7-2b_8=0,$

$2b_5+2b_6+2b_7-2b_8=0.$ 

Whence it follows that
$\dim\,Z_{0}^3({\mathfrak{g}})=3.$ By~(3),
$$\dim\,H_0^3({\mathfrak{g}})=\dim\,\overline{Z}_0^3({\mathfrak{g}})+
\dim\,\overline{Z}_0^2({\mathfrak{g}})-\dim\,\overline{C}_0^2({\mathfrak{g}})=3+2-4=1.$$
Therefore $[H_0^3({\mathfrak{g}}):k]=1.$

The weight subspaces $\overline {C}^3_{3\lambda_1}({\mathfrak{g}}),$
$\overline{C}^3_{3\lambda_2}({\mathfrak{g}})$ are two-dimensional and
span respectively with $3$-cochains: $$h_1^*\wedge f_1^*\wedge
f_3^*,\,h_2^*\wedge f_1^*\wedge f_3^*,\, h_1^*\wedge
f_2^*\wedge f_3^*,\, h_2^*\wedge f_2^*\wedge f_3^*.$$ Using the
cocycle condition, we get
$\dim\,\overline{Z}^3_{3\lambda_1}({\mathfrak{g}})=\dim\,
\overline{Z}^3_{3\lambda_2}({\mathfrak{g}})=1.$ So,
$H^3({\mathfrak{g}})\cong L(1,0)^{(1)}\oplus L(0,1)^{(1)}\oplus k.$

Now we prove that $H^4({\mathfrak{g}})=0.$ It's obvious that
$\prod(\overline{C}^4({\mathfrak{g}}))=\prod(\overline{C}^3({\mathfrak{g}})).$
Therefore we consider only the weight subspaces
$$\overline{C}^4_{0}({\mathfrak{g}}),\,\overline{C}^4_{3\omega_1}({\mathfrak{g}}),\,\overline{C}^4_{3\omega_2}({\mathfrak{g}}).$$

\medskip
\noindent The subspace $\overline{C}^4_{0}({\mathfrak{g}})$ is
$10$-dimensional and spans by the cochains
$$h_1^*\wedge h_2^*\wedge e_1^*\wedge f_1^*,\,h_1^*\wedge h_2^*\wedge e_2^*\wedge f_2^*,\,
h_1^*\wedge h_2^*\wedge e_3^*\wedge f_3^*,\, h_1^*\wedge e_1^*\wedge
e_2^*\wedge f_3^*,$$
$$h_2^*\wedge e_1^*\wedge e_2^*\wedge f_3^*,\,h_1^*\wedge e_3^*\wedge f_1^*\wedge f_2^*,\,h_2^*\wedge e_3^*\wedge f_1^*\wedge f_2^*,$$
$$e_1^*\wedge e_2^*\wedge f_1^*\wedge f_2^*,\, e_1^*\wedge e_3^*\wedge f_1^*\wedge f_3^*,\,e_2^*\wedge e_3^*\wedge f_2^*\wedge f_3^*.$$
Suppose that the linear combination of these vectors coefficients
$b_i,\,i=1,\cdots,10$ respectively is a $4$-cocycle. Then
$b_1=b_2=b_3=0,$ $b_4=b_6,$ $b_5=b_7.$ Whence it follows that
$\dim\,\overline{z}^4_{0}({\mathfrak{g}})=5.$ By~(3),
$\dim\,H_{0}^4({\mathfrak{g}})=5+3-8=0.$ Therefore,
$[H^4({\mathfrak{g}}):k]=0.$

It's obvious that
$\overline{Z}^4_{3\omega_1}({\mathfrak{g}},M))=\overline{Z}^4_{3\omega_1}({\mathfrak{g}},M))=1.$
Then by~(3),
$$\overline{H}^4_{3\omega_1}({\mathfrak{g}},M))=\overline{H}^4_{3\omega_1}({\mathfrak{g}},M))=1+1-2=0.$$
So, $H^4({\mathfrak{g}})=0.$

Using~(2) and the statements $(b),\,(c),$ we get the statements
$(e),\,(f)$ respectively. The proof of Lemma~\ref{l1} is complete.

Now let $M=L(1,1).$

\begin{lm}\label{l2} There are the following isomorphisms of
$SL_3(k)$-modules:

(a) $H^1({\mathfrak{g}},L(1,1))\cong L(1,0)^{(1)}\oplus
L(0,1)^{(1)}\oplus k;$

(b) $H^3({\mathfrak{g}},L(1,1))\cong H^0(1,1)^{(1)};$

(c) $H^4({\mathfrak{g}},L(1,1))\cong 2H^0(1,1)^{(1)};$

(d) $H^5({\mathfrak{g}},L(1,1))\cong H^0(1,1)^{(1)};$

(e) $H^7({\mathfrak{g}},L(1,1))\cong L(1,0)^{(1)}\oplus
L(0,1)^{(1)}\oplus k.$

In other cases $H^n({\mathfrak{g}},L(1,1))=0.$
\end{lm}

{\bf Proof.} The calculations similar to the previous Lemma~\ref{l1}
yield:

1)
$\prod(\overline{C}^0({\mathfrak{g}},L(1,1)))=\prod(\overline{C}^8({\mathfrak{g}},L(1,1)))=\{0\},$

\noindent
$\prod(\overline{C}^i({\mathfrak{g}},L(1,1)))=
\{0,\pm3\omega_1,\pm3(\omega_1-\omega_2),\pm3\omega_2\}$
for $i=1,2,6,7,$

\noindent
$\prod(\overline{C}^j({\mathfrak{g}},L(1,1)))=\prod(\overline{C}^1({\mathfrak{g}},L(1,1)))\cup$

\noindent
$\{\pm3(\omega_1+\omega_2),\pm3(2\omega_1-\omega_2),\pm3(-\omega_1+2\omega_2)\}$
for $j=3,4;$

2)
$\dim\,\overline{C}^0_0({\mathfrak{g}},L(1,1))=\dim\,\overline{C}_0^8({\mathfrak{g}},L(1,1))=1,$
$\dim\,\overline{C}_0^1({\mathfrak{g}},L(1,1))=\dim\,\overline{C}_0^7({\mathfrak{g}},L(1,1))=8,$

\noindent
$\dim\,\overline{C}_0^2({\mathfrak{g}},L(1,1))=\dim\,\overline{C}_0^6({\mathfrak{g}},L(1,1))=22,$

\noindent
$\dim\,\overline{C}_0^3({\mathfrak{g}},L(1,1))=\dim\,\overline{C}_0^5({\mathfrak{g}},L(1,1))=38,$

\noindent $\dim\,\overline{C}_0^4({\mathfrak{g}},L(1,1))=44;$

3)
$\dim\,\overline{C}^0_{3\omega_i}({\mathfrak{g}},L(1,1))=\dim\,\overline{C}_{3\omega_i}^8({\mathfrak{g}},L(1,1))=0,$
$\dim\,\overline{C}_{3\omega_i}^1({\mathfrak{g}},L(1,1))=\dim\,\overline{C}_{3\omega_i}^7({\mathfrak{g}},L(1,1))=2,$

\noindent
$\dim\,\overline{C}_{3\omega_i}^2({\mathfrak{g}},L(1,1))=\dim\,\overline{C}_{3\omega_i}^6({\mathfrak{g}},L(1,1))=7,$
$\dim\,\overline{C}_{3\omega_i}^3({\mathfrak{g}},L(1,1))=\dim\,\overline{C}_{3\omega_i}^5({\mathfrak{g}},L(1,1))=14,$

\noindent $\dim\,\overline{C}_{3\omega_i}^4({\mathfrak{g}},L(1,1))=18$
for $i=1,2;$

4)
$\dim\,\overline{C}^0_{3(\omega_1+\omega_2)}({\mathfrak{g}},L(1,1))=\dim\,\overline{C}_{3(\omega_1+\omega_2)}^8({\mathfrak{g}},L(1,1))=0,$

\noindent
$\dim\,\overline{C}_{3(\omega_1+\omega_2)}^1({\mathfrak{g}},L(1,1))=\dim\,\overline{C}_{3(\omega_1+\omega_2)}^7({\mathfrak{g}},L(1,1))=0,$

\noindent
$\dim\,\overline{C}_{3(\omega_1+\omega_2)}^2({\mathfrak{g}},L(1,1))=\dim\,\overline{C}_{3(\omega_1+\omega_2)}^6({\mathfrak{g}},L(1,1))=0,$

\noindent
$\dim\,\overline{C}_{3(\omega_1+\omega_2)}^3({\mathfrak{g}},L(1,1))=\dim\,\overline{C}_{3(\omega_1+\omega_2)}^5({\mathfrak{g}},L(1,1))=1,$

\noindent
$\dim\,\overline{C}_{3(\omega_1+\omega_2)}^4({\mathfrak{g}},L(1,1))=2;$

5)
$\dim\,\overline{Z}^0_0({\mathfrak{g}},L(1,1))=\dim\,\overline{C}_0^8({\mathfrak{g}},L(1,1))=0,$
$\dim\,\overline{Z}_0^1({\mathfrak{g}},L(1,1))=\dim\,\overline{Z}_0^7({\mathfrak{g}},L(1,1))=2,$

\noindent
$\dim\,\overline{Z}_0^2({\mathfrak{g}},L(1,1))=\dim\,\overline{Z}_0^6({\mathfrak{g}},L(1,1))=6,$

\noindent
$\dim\,\overline{Z}_0^3({\mathfrak{g}},L(1,1))=\dim\,\overline{Z}_0^5({\mathfrak{g}},L(1,1))=18,$

\noindent $\dim\,\overline{Z}_0^4({\mathfrak{g}},L(1,1))=24;$

6)
$\dim\,\overline{Z}^0_{3\omega_i}({\mathfrak{g}},L(1,1))=\dim\,\overline{Z}_{3\omega_i}^8({\mathfrak{g}},L(1,1))=0,$
$\dim\,\overline{Z}_{3\omega_i}^1({\mathfrak{g}},L(1,1))=\dim\,\overline{Z}_{3\omega_i}^7({\mathfrak{g}},L(1,1))=1,$

\noindent
$\dim\,\overline{Z}_{3\omega_i}^2({\mathfrak{g}},L(1,1))=\dim\,\overline{Z}_{3\omega_i}^6({\mathfrak{g}},L(1,1))=1,$
$\dim\,\overline{C}_{3\omega_i}^3({\mathfrak{g}},L(1,1))=\dim\,\overline{C}_{3\omega_i}^5({\mathfrak{g}},L(1,1))=6,$

\noindent $\dim\,\overline{C}_{3\omega_i}^4({\mathfrak{g}},L(1,1))=8$
for $i=1,2;$

7)
$\dim\,\overline{Z}^0_{3(\omega_1+\omega_2)}({\mathfrak{g}},L(1,1))=\dim\,\overline{Z}_{3(\omega_1+\omega_2)}^8({\mathfrak{g}},L(1,1))=0,$

\noindent
$\dim\,\overline{Z}_{3(\omega_1+\omega_2)}^1({\mathfrak{g}},L(1,1))=\dim\,\overline{Z}_{3(\omega_1+\omega_2)}^7({\mathfrak{g}},L(1,1))=0,$

\noindent
$\dim\,\overline{Z}_{3(\omega_1+\omega_2)}^2({\mathfrak{g}},L(1,1))=\dim\,\overline{Z}_{3(\omega_1+\omega_2)}^6({\mathfrak{g}},L(1,1))=0,$

\noindent
$\dim\,\overline{Z}_{3(\omega_1+\omega_2)}^3({\mathfrak{g}},L(1,1))=\dim\,\overline{Z}_{3(\omega_1+\omega_2)}^5({\mathfrak{g}},L(1,1))=1,$

\noindent
$\dim\,\overline{Z}_{3(\omega_1+\omega_2)}^4({\mathfrak{g}},L(1,1))=2.$

Then, by~(3), $\dim\,\overline{H}^n_{\mu}({\mathfrak{g}},L(1,1))=0$
except in the following cases:

i)
$\dim\,\overline{H}^1_0({\mathfrak{g}},L(1,1))=\dim\,\overline{H}_0^7({\mathfrak{g}},L(1,1))=1,$
$\dim\,\overline{H}^3_0({\mathfrak{g}},L(1,1))=\dim\,\overline{H}_0^5({\mathfrak{g}},L(1,1))=2,$

\noindent $\dim\,\overline{H}^4_0({\mathfrak{g}},L(1,1))=4;$

ii)
$\dim\,\overline{H}^1_{3\omega_i}({\mathfrak{g}},L(1,1))=\dim\,\overline{H}_{3\omega_i}^7({\mathfrak{g}},L(1,1))=1$
for $i=1,2;$

iii)
$\dim\,\overline{H}^3_{3(\omega_1+\omega_2)}({\mathfrak{g}},L(1,1))=\dim\,\overline{H}_{3(\omega_1+\omega_2)}^5({\mathfrak{g}},L(1,1))=1,$

\noindent
$\dim\,\overline{H}^4_{3(\omega_1+\omega_2)}({\mathfrak{g}},L(1,1))=2.$

\noindent Analyzing the dimensions of the weight subspaces of the
corresponding cohomology groups, we obtain the required statements
of Lemma~\ref{l2}. The proof of Lemma~\ref{l2} is complete.

Combining the results of Lemmas~\ref{l1} and~\ref{l2}, we obtain all the
statements of Theorem~\ref{th1}.

\section{Cohomlogy of the adjoint module}\label{s3}

Using Theorem~\ref{th1}, we can easily compute the cohomology of the adjoint
module for ${\mathfrak{g}}.$ There is the following short exact
sequence of ${\mathfrak{g}}$-modules:
$$0\rightarrow k\rightarrow {\mathfrak{g}}\rightarrow L(1,1)\rightarrow
0.$$ Consider the corresponding long exact cohomological sequence of
$SL_3(k)$-modules
$$\cdots \rightarrow H^{n-1}({\mathfrak{g}},L(1,1))\rightarrow H^n({\mathfrak{g}})\rightarrow H^n({\mathfrak{g}},{\mathfrak{g}})
\rightarrow $$$$H^n({\mathfrak{g}},L(1,1))\rightarrow
H^{n+1}({\mathfrak{g}})\rightarrow\cdots$$

It is known that $H^2({\mathfrak{g}},{\mathfrak{g}})=0$ \cite{14}. Then,
according to Theorem~\ref{th1}, the last long exact cohomological sequence
splits into the following five exact sequences:
$$0 \rightarrow H^0({\mathfrak{g}})\rightarrow H^0({\mathfrak{g}},{\mathfrak{g}})
\rightarrow 0,$$ $$0\rightarrow H^1({\mathfrak{g}},{\mathfrak{g}})
\rightarrow H^1({\mathfrak{g}},L(1,1))\rightarrow
H^{2}({\mathfrak{g}})\rightarrow 0,$$
$$0\rightarrow H^3({\mathfrak{g}})\rightarrow
H^3({\mathfrak{g}},{\mathfrak{g}}) \rightarrow
H^3({\mathfrak{g}},L(1,1))\rightarrow 0,$$
$$0\rightarrow H^4({\mathfrak{g}},{\mathfrak{g}})
\rightarrow H^4({\mathfrak{g}},L(1,1))\rightarrow
H^{5}({\mathfrak{g}})\rightarrow H^5({\mathfrak{g}},{\mathfrak{g}})
\rightarrow $$$$H^5({\mathfrak{g}},L(1,1))\rightarrow
H^{6}({\mathfrak{g}})\rightarrow H^6({\mathfrak{g}},{\mathfrak{g}})
\rightarrow 0,$$
$$0\rightarrow H^7({\mathfrak{g}},{\mathfrak{g}})
\rightarrow H^7({\mathfrak{g}},L(1,1))\rightarrow
H^{8}({\mathfrak{g}})\rightarrow H^8({\mathfrak{g}},{\mathfrak{g}})
\rightarrow 0.$$

The first three short exact sequences yield the following
isomorphisms of $SL_3(k)$-modules respectively:
$$H^0({\mathfrak{g}},{\mathfrak{g}})\cong k,\,H^1({\mathfrak{g}},{\mathfrak{g}})\cong k,$$$$H^3({\mathfrak{g}},{\mathfrak{g}})\cong L(1,0)^{(1)}\oplus
L(0,1)^{(1)}\oplus H^0(1,1)^{(1)}\oplus k.
$$
Since $3(\omega_1+\omega_2)\notin \prod(H^i({\mathfrak{g}}))$ for
$i=5,6,$ then the fourth exact sequence splits and yields the
following isomorphisms:
$$H^4({\mathfrak{g}},{\mathfrak{g}})\cong H^4({\mathfrak{g}},L(1,1))\cong H^0(1,1)^{(1)},$$$$H^5({\mathfrak{g}},{\mathfrak{g}})\cong L(1,0)^{(1)}\oplus
L(0,1)^{(1)}\oplus H^0(1,1)^{(1)}\oplus k,$$
$$H^6({\mathfrak{g}},{\mathfrak{g}})\cong L(1,0)^{(1)}\oplus
L(0,1)^{(1)}.$$ Similarly to the previous case, from the last exact
sequence we obtain $$H^7({\mathfrak{g}},{\mathfrak{g}})\cong
L(1,0)^{(1)}\oplus L(0,1)^{(1)}\oplus k,
H^8({\mathfrak{g}},{\mathfrak{g}})\cong k.$$

Thus, we get the following

\begin{prp}\label{p1} Let ${\mathfrak{g}}$ be a classical Lie algebra
of type $A_2$ over an algebraically closed field $k$ of
characteristic $p=3.$ Then there are the following isomorphisms of
$SL_3(k)$-modules:

(a) $H^0({\mathfrak{g}},{\mathfrak{g}})\cong
H^1({\mathfrak{g}},{\mathfrak{g}})\cong
H^8({\mathfrak{g}},{\mathfrak{g}})\cong k;$

(b) $H^3({\mathfrak{g}},{\mathfrak{g}})\cong
H^5({\mathfrak{g}},{\mathfrak{g}})\cong L(1,0)^{(1)}\oplus
L(0,1)^{(1)}\oplus H^0(1,1)^{(1)}\oplus k;$

(c) $H^4({\mathfrak{g}},{\mathfrak{g}})\cong 2H^0(1,1)^{(1)};$

(d) $H^6({\mathfrak{g}},{\mathfrak{g}})\cong
 L(1,0)^{(1)}\oplus L(0,1)^{(1)};$

(e) $H^7({\mathfrak{g}},{\mathfrak{g}})\cong L(1,0)^{(1)}\oplus
L(0,1)^{(1)}\oplus k.$

In other cases $H^n({\mathfrak{g}},{\mathfrak{g}})=0.$
\end{prp}

\section{ Cohomlogy for $\overline{A_2}$}\label{s4}

Recall that $\overline{A_2}$ is the quotient algebra of the
classical Lie algebra of type $A_2$ over an algebraically closed
field of characteristic $p=3$ by the center. In this section we
compute cohomology of the simple Lie algebra $\overline{A_2}$ with
coefficients in the simple modules.

First, we consider an arbitrary Lie algebra $\mathfrak{g}$ with the
center $C_{\mathfrak{g}}$ such that the corresponding quotient
algebra is a simple algebra. The following result will immediately
lead to our goal.

\begin{lm}\label{l3} Let $\overline{\mathfrak{g}}$ be a simple quotient
Lie algebra of a Lie algebra  $\mathfrak{g}$ by the center
$C_{\mathfrak{g}}.$ Then
$H^n(\overline{\mathfrak{g}},\overline{\mathfrak{g}})\cong
H^n(\mathfrak{g},\overline{\mathfrak{g}})$ for all $n>0.$
\end{lm}

{\bf Proof.} The space $\overline{\mathfrak{g}}$ can be equipped
with the structure  of a module over each of the Lie algebras
$C_{\mathfrak{g}},$ $\mathfrak{g}$ and $\overline{\mathfrak{g}}:$

$C_{\mathfrak{g}}\times \overline{\mathfrak{g}}\rightarrow
\overline{\mathfrak{g}}, \;(c,\overline{a})\mapsto
\mu(c)\overline{a},$ where $\mu$ is a nonzero linear form on
$C_{\mathfrak{g}};$

$\mathfrak{g}\times \overline{\mathfrak{g}}\rightarrow
\overline{\mathfrak{g}}, \;(a_1,\overline{a}_2)\mapsto
\overline{[a_1,a_2]},\;a_1\in
\mathfrak{g},\;\overline{a}_2\in\overline{\mathfrak{g}};$

$\overline{\mathfrak{g}}\times \overline{\mathfrak{g}}\rightarrow
\overline{\mathfrak{g}}, \;(\overline{a}_1,\overline{a}_2)\mapsto
\overline{[a_1,a_2]},\;\overline{a}_1,\;\overline{a}_2\in\overline{\mathfrak{g}}.$

The short exact sequence of cochain complexes
$$0\rightarrow (C^*(C_{\mathfrak{g}},\overline{\mathfrak{g}}),d)\rightarrow
(C^*({\mathfrak{g}},\overline{\mathfrak{g}}),d)\rightarrow
(C^*(\overline{\mathfrak{g}},\overline{\mathfrak{g}}),d)\rightarrow
0$$ gives a long exact cohomological sequence
$$\cdots \rightarrow
H^{n-1}(C_{\mathfrak{g}},\overline{\mathfrak{g}})\rightarrow
H^n(\overline{\mathfrak{g}},\overline{\mathfrak{g}})\rightarrow
H^n({\mathfrak{g}},\overline{\mathfrak{g}})\rightarrow
H^n(C_{\mathfrak{g}},\overline{\mathfrak{g}})\rightarrow\cdots.
$$
Since $H^n(C_{\mathfrak{g}},\overline{\mathfrak{g}})=0$ for all
$n\geq 0$ \cite[כוללא~4.2]{15}, it follows from the fact that last
cohomological sequence is exact that
$H^n(\overline{\mathfrak{g}},\overline{\mathfrak{g}})\cong
H^n(\mathfrak{g},\overline{\mathfrak{g}})$ for all $n>0.$ The proof
of Lemma 3 is complete.

\begin{rmk}\label{r1} A special case of Lemma~3 for $n=1$ was proved
in~\cite{7}. Using Lemma~\ref{l3} to Theorem~\ref{th1}, we obtain a complete description
of the cohomology of a simple Lie algebra $\overline{A_2}$ with
coefficients in simple modules.
\end{rmk}

\bigskip
\textsf{This research was funded by Science Committee grant of the
Ministry of Education and Science of the Republic of Kazakhstan
(Grant No. AP08855935).}

\end{document}